\documentclass{article}
\usepackage{latexsym}
\usepackage{amssymb}
\usepackage{amsfonts}                     
\usepackage{amsmath}
\usepackage[all]{xy}
\usepackage{color}
\usepackage{mathrsfs}
\usepackage[utf8]{inputenc}
\usepackage[T1]{fontenc}
\usepackage{amscd}
\usepackage{amsxtra} 
\newtheorem{thm}{\bf Theorem}[section]
\newtheorem{cor}[thm]{\bf Corollary}
\newtheorem{lem}[thm]{\bf Lemma}
\newtheorem{prop}[thm]{\bf Proposition}
\newtheorem{defn}[thm]{\bf Definition}

\newtheorem{rems}[thm]{\bf Remarks}
\newtheorem{exmp}[thm]{\bf Example}

\newcommand{\field}[1]{\mathbb{#1}}

\newcommand{\Q }{\field{Q}}
\newcommand{\Z }{\field{Z}}

\def\C{{\cal C}}
\def\E{{\cal E}}
\def\F{{\cal F}}

\def\X{{\cal X}}
\def\Y{{\cal Y}}

\def\Ext{{\rm Ext}}
\def\Hom{{\rm Hom}}
\def\Im{{\rm Im}}

\def\proof{{\parindent0pt {\bf Proof.\ }}}
\def\Tor{{\rm Tor}}
\def\coker{{\rm coker}}

\def\Card{{\rm Card}}

\newcommand{\cqfd}
{\hspace{1cm}
\rule{2mm}{2mm}%
\medbreak%
\par%
}
\begin{document}
\title{$S$-flat cotorsion pair }
\author{Driss Bennis and Ayoub Bouziri}

\date{}
 
\maketitle
\begin{abstract} 
 
Let $R$ be a commutative ring, and let $S$ be a multiplicative subset of $R$. In this paper, we investigate the notion of $S$-cotorsion modules.   An $R$-module $C$ is called $S$-cotorsion if $\Ext^{1}_{R}(F,C) = 0$ for every $S$-flat $R$-module $F$. Among other results, we establish that the pair $(S\F, S\C)$, where $S\F$ denotes the class of all $S$-flat $R$-modules and $S\C$ denotes the class of all $S$-cotorsion modules, forms a hereditary perfect cotorsion pair. As applications, we provide characterizations of $S$-perfect rings in terms of $S$-cotorsion modules. We conclude the paper with results on $S\F$-preenvelopes. Namely, we prove that if every module has an $S\F$-preenvelope, then $R$ is $S$-coherent. Furthermore, we establish the converse under the condition that $R_S$ is a finitely presented $R$-module.
\end{abstract}

\medskip



\section{Introduction}
\hskip .5cm  Throughout this paper, $R$ is a commutative ring with identity, all modules are unitary and $S$ is a multiplicative subset of $R$; that is,  $1 \in S$ and $s_1s_2 \in S$ for any $s_1,s_2 \in S$. Unless explicitly stated otherwise, when we consider a multiplicative subset $S$ of $R$, we implicitly suppose that $0\notin S$. This will be used in the sequel without explicit mention. 

\medskip
 Let $\X$ be a class of $R$-modules and $M$ an $R$-module. Following Enochs \cite{Eno1}, we say that a homomorphism $ \phi : M\to X$  is an $\X$-preenvelope if $X \in \X$ and the abelian group homomorphism $ \Hom_R(\phi, X') : \Hom_R(X, X')\to\Hom_R(M, X')$ is an epimorphism for each $X' \in \X$. An $\X$-preenvelope $\phi : M \to X$ is said to be an $\X$-envelope if every endomorphism $g: X \to X$ such that $g\phi = \phi$ is an isomorphism.  We will denote by \begin{center} $\X ^\bot = \lbrace X : \Ext^1_R(Y, X) = 0$ for all Y$ \in\X \rbrace$\\
$^\bot\X = \lbrace X : \Ext^1_R(X,Y) = 0$ for all Y$ \in \X\rbrace$
\end{center} the right orthogonal class and the left orthogonal class of $\X$, respectively.

\medskip 
Following \cite{Eno2}, an epimorphism  $ \alpha : M\to X$ with $X \in \X$ is said to be a special $\X$-preenvelope of $M$ if $ \coker(\alpha)\in \,^\bot\X$.  Dually, we have the definitions of a (special) $\X$-precover and an $\X$-cover. $\X$-envelopes ($\X$-covers) may not exist in general, but if they exist, they are unique up to isomorphism.

\medskip
A pair $(\X, \Y)$ of classes of $R$-modules is called a cotorsion pair (or cotorsion theory \cite{Eno2}) if $\X ^\bot =\Y$ and $^\bot \Y =\X $. A cotorsion pair $(\X, \Y)$ is called perfect if every $R$-module has an $\X$-cover and a $\Y$-envelope \cite{Eno3}. A cotorsion pair $(\X, \Y)$ is called complete  \cite{Trl1} if for any $R$-module $M$, there are exact sequences $0\to M  \to Y \to X \to 0$ with $ Y\in \Y $ and $X \in \X$, and $0 \to Y' \to X'\to \ M \to 0$ with $ X'\in \X$ and $Y'\in \Y$. A perfect cotorsion pair is always complete by Wakamatsu's Lemmas \cite[Section 2.1]{Xu1}.  A cotorsion pair $(\X, \Y)$ is said to be hereditary if whenever  $0\to X'  \to X \to X'' \to 0$ is exact with $ X, X''\in \X$, then $X'$ is also in $\X$ \cite{Eno3}.  According to \cite[Proposition 1.2]{Eno3}, a cotorsion pair $(\X, \Y)$ is hereditary if and only if, whenever $0\to Y'  \to Y \to Y'' \to 0$ is exact with $ Y, Y'\in \Y$, $Y''$ is also in $\Y$.

\medskip
In \cite{Ben1}, Bennis and El Hajoui investigated an $S$-version of finitely presented modules and coherent rings which are called, respectively, $S$-finitely presented modules and $S$-coherent rings. An $R$-module $M$ is said to be $S$-finitely presented if there exists an exact sequence of $R$-modules $0 \to K\to L\to  M \to 0$, where $L$ is a  finitely generated free $R$-module and $K$ is an $S$-finite $R$-module; that is, there exist a finitely generated submodule $F$ of $K$ and $s \in S$ such that $sK \subseteq F$. Moreover, a  commutative ring $R$ is called $S$-coherent if every finitely generated ideal of $R$ is $S$-finitely presented. They showed that the $S$-coherent rings have a similar characterization to the classical one given by Chase for coherent rings  \cite[Theorem 3.8]{Cha1}. Subsequently, they asked whether there exists an $S$-version of Chase's theorem \cite[Theorem 2.1]{Cha1}. In other words, how to define an $S$-version of flatness that characterizes $S$-coherent rings similarly to the classical case?  This problem was solved by the  notion  of $S$-flat modules in \cite{Qi1}. 

\medskip
  Recently, we have introduced and studied the notion of $S$-perfect rings. A ring $R$ is said to be $S$-perfect if any $S$-flat $R$-module is projective \cite{Bou1}. Several characterizations of $S$-perfect rings are given in \cite{Bou1}. In this work, we aim to contribute new characterizations in terms of $S$-cotorsion modules (see Proposition \ref{3-prop-S-perfect} and Theorem \ref{3-thm-S-perfect}).

\medskip
The organization of the paper is as follows: In Section 2, several elementary properties of $S$-flat modules are obtained. The concept of $S$-cotorsion modules, which is different from $S$-cotorsion modules in the sense of \cite{Ass1}, is first introduced in Section $3$. An $R$-module $M$ is said to be $S$-cotorsion if $\Ext^1_R(F, M) = 0$ for any $S$-flat $R$-module $F$. We prove that the pair ($S\F,S\C)$, where $S\F$ is the class of all $S$-flat $R$-modules and $S\C$ is the class of all $S$-cotorsion modules, is a hereditary perfect cotorsion pair (see Theorem \ref{3-thm-(SF,SC)-perfect-cotorsion-pair}). We show that $S$-perfect rings are characterized in terms of $S$-cotorsion module (see Proposition \ref{3-prop-S-perfect} and Theorem \ref{3-thm-S-perfect}). Other results, on $S$-cotorsion envelopes, represent the $S$-counterpart of that of the cotorsion envelopes \cite{Mao1}. In Section 4, we deal with $S\F$-preenvelope, we prove that if any module has an $S\F$-preenvelope, then $R$ is $S$-coherent (see Corollary \ref{4-cor-S-coherent}), and we prove the converse when $R_S$ is a finitely presented $R$-module (see Proposition \ref{4-prop-S-coherent-Rs-EP}). 

\medskip
From now on, we will write $S\F$ for the class of all $S$-flat $R$-modules and $S\C$ for the class of all $S$-cotorsion $R$-modules. For an $R$-module $M$, we write $M_S$ to indicate the localization of $M$ at $S$. $\epsilon_M : S\F(M)\to M$ and $ \sigma_M : M\to S\C(M)$ will denote an $S\F$-cover and an $S\C$-envelope of $M$, respectively. Sometimes we just call $S\C(M)$ an $S\C$-envelope of $M$. We use $N \leq M$ and $N \leq_e M$ to mean that $N$ is a submodule and an essential submodule of $M$, respectively. Finally, the character module $\Hom_\Z(M,\Q/\Z)$ will be denoted by $M^+$. 


\section{S-flat modules}
 Recall that an $R$-module $M$ is said to be  $S$-flat if, for any finitely generated ideal $I$ of $R$, the induced homomorphism $(I \otimes_{R} M)_S \to (R \otimes_{R} M)_S $ is a monomorphism; equivalently, $M_{S}\cong R_S\otimes_R M$ is a flat $R_{S}$-module \cite[Proposition 2.6]{Qi1}. 
 
 It is  well-know that the class $\F$ of all flat modules is closed under pure submodules, pure quotient modules, extensions and direct limits. Here we have the corresponding result for the class of all $S$-flat modules $S\F.$

\begin{lem}\label{2-lem-S-flat-class-properties}
$S\F$ is closed under extensions, direct sums, direct summands, direct limits, pure submodules and pure quotients.
\end{lem}

\proof If $0\to A\to B \to C \to 0 $ is a (pure) exact sequence of $R$-modules, then the induced sequence $0\to A_{S}\to B_{S} \to C_{S} \to 0 $ is a (pure) exact sequence of $R_{S}$-modules. Thus, all properties follow from their validity for the class of flat $R_S$-modules and the fact that they are preserved by the functor $R_S\otimes_R(-)$.  \cqfd


Recall from \cite[Definition 1.6.10]{Wan1} the following definition: 

\begin{defn}[\cite{Wan1}, Definition 1.6.10]\label{2-def-S-torsion} Let $M$ be an $R$-module. Set
$$tor_S(M) = \{x \in M | \text{ there exists } s\in S \text{ such that } sx = 0\}.$$
Then $tor_S(M)$ is a submodule of M, called the $S$-torsion submodule of $M$ and $M$ is
called an $S$-torsion if $tor_S(M) = M$. 
\end{defn}

 The canonical ring homomorphism $\theta: R \to R_{S}$ makes every $R_{S}$-module an $R$-module via the formula $r.m = \frac{r}{1}.m$, where $r\in R $ and $m\in M$. Recall from \cite{Rot1} the following lemma that we frequently use in this paper.

\begin{lem}[\cite{Rot1}, Corollary 4.79]\label{2-lem-localiz-isomorphism}
Every $R_S$-module $M$ is naturally isomorphic to its localization $M_S$ as $R_S$-modules
\end{lem}

Recall that a sequence $0 \to A \to B\to C \to 0$ is $S$-exact if the induced sequence
$0 \to A_S \to A_S \to C_S \to 0$ is exact \cite[Definition 2.1]{Qi1}. Since $R_S$ is a flat $R$-module \cite[Theorem 4.80]{Rot1}, every exact sequence is $S$-exact. The following lemma follows form the standard arguments: 
\begin{lem}\label{2-lem-S-flat-exact-chara} The following assertions are equivalent for an $R$-module $M$:
\begin{enumerate}
\item $M$ is an  $S$-flat $R$-module.
\item For every  $S$-exact sequence $0\to A\to B\to C\to 0$ of $R$-modules, the induced sequence  $0\to A\otimes_{R} M\to B\otimes_{R} M\to C\otimes_{R} M\to 0$ is  $S$-exact.
\item For every short exact sequence $0\to A\to B\to C\to 0$ of $R$-modules, the induced sequence  $0\to A\otimes_{R} M\to B\otimes_{R} M\to C\otimes_{R} M\to 0$ is  $S$-exact.
\item For every  $S$-exact sequence $0\to K\to L\to M\to 0$ of $R$-modules, the induced sequence $0\to K_{S}\to L_{S}\to M_{S}\to 0$ is a pure-exact sequence of $R_{S}$-modules.
\item For every exact sequence $0\to K\to L\to M\to 0$ of $R$-module, the induced sequence $0\to K_{S}\to L_{S}\to M_{S}\to 0$ is pure-exact sequence of $R_{S}$-modules.
\item $ \Tor_{R}^{n}(M,N)$ is $S$-torsion for any $R$-module $N$ and $n\geq 1 $.
\end{enumerate}
\end{lem}

It is well-known that flat $R$-modules can be characterized in terms of $\Tor$ functor \cite[Theorem 1.2.1]{Gla1}. Now, we explore similar properties of $S$-flat $R$-modules in relation to the $\Tor$.

\begin{prop}\label{2-prop-S-flat-tor-chara} The following assertions are equivalent for an $R$-module $M$:
\begin{enumerate}
\item $M$ is  $S$-flat;
\item $ \Tor^1_{R}(M,N)=0$ for any $R_{S}$-module $N$;
\item $ \Tor^1_{R}(M,N_{S})=0$ for any $R$-module $N$;
\item $ \Tor^n_{R}(M,N)=0$ for any $R_{S}$-module $N$ and $n\geq 1$;
\item $ \Tor^n_{R}(M,N_{S})=0$ for any $R$-module $N$ and $n\geq 1$.
\item $ \Tor^n_{R}(M,(R/I)_{S})=0$ for any (finitely generated) ideal $I$ of $R$ and $n\geq 1$.
\item $ \Tor^1_{R}(M,(R/I)_{S})=0$ for any (finitely generated) ideal $I$ of $R$.

\end{enumerate}
\end{prop}

\proof  $ 1.\Rightarrow 2. $ Let $N$ be an $R_{S}$-module, and let $ 0\to K\to P\to N\to 0$ be an exact sequence of $R$-modules, where $P$ is a projective $R$-module. Notice that $N_{S}\cong N$.  We have the exact sequence of $R_{S}$-modules
$$ 0\to K_S\to P_S\to N\to 0,$$
which yields the exactness of the sequence
$$ \Tor^{1}_{R_{S}}(M_{S},N)=0\to M_S\otimes_{R_{S}} K_{S}\to M_S\otimes_{R_{S}} P_{S}\to M_S\otimes_{R_{S}} N\to 0$$
which gives rise to the exactness of the sequence
$$  0\to M\otimes_{R} K_{S}\to M\otimes_{R} P_{S}\to M\otimes_{R} N\to 0. $$
On the other hand, the sequence 
$$  0\to\Tor^{1}_{R}(M,N)\to M\otimes_{R} K_{S}\to M\otimes_{R} P_{S}\to M\otimes_{R} N\to 0 $$
is exact. Thus, $\Tor^{1}_{R}(M,N)=0$, as desired.

$ 2.\Leftrightarrow 3.$ and $ 4.\Leftrightarrow 5.$ are clear.

$ 3.\Rightarrow 5.$ Let $N$ be an $R$-module. The proof is by induction on $n \geq 1$. There is an exact sequence $0 \to K_{S} \to P_{S} \to N_{S} \to 0$ with $P$ a free $R$-module. For the inductive step, we use the long exact sequence theorem to obtain the exactness of 
$$0 =\Tor_R^{ n+1}(M, P_S)\to\Tor_R^{n+1}(M, N_S) \to\Tor_R^n (M, K_S) \to\Tor_R^n(M,P_S) = 0.$$
But, by induction, $\Tor_R^n (M, K_S)=0$, then $\Tor_R^{n+1}(M, N_S)=0$.

$ 5.\Rightarrow 6.$ and $ 6.\Rightarrow 7.$ are clear.

$7.\Rightarrow 1.$ Let $I$ be a finitely generated ideal of $R_{S}$. We can set $I= J_S$, where $J$ is a finitely generated ideal of $R$. We  have  $R_S/I\cong (R/J)_S$. Then
$$ \Tor_{R_{S}}^{1}(M_{S},R_S/I)\cong\Tor_{R}^{1}(M, R/J )_{S} $$ \cite[Proposition 7.17]{Rot1}. By $(7)$  the right hand is zero. Thus, $M_{S}$ is a flat $R_{S}$-module. Then $M$ is an $S$-flat $R$-module \cite[Proposition 2.6]{Qi1}. \cqfd

\begin{cor}\label{2-cor-S-flat-Rs-module} The following are equivalent for an $R_S$-module $M$:
\begin{enumerate}
\item $M$ is an $S$-flat $R$-module.
\item $M$ is a flat $R_S$-module.
\item $M$ is a flat $R$-module.
\end{enumerate} 
\end{cor}

\begin{cor}\label{2-cor-of-S-flat-tor-chara-kernel-extension}
Let $0\to K\to L\to M\to 0$ be an exact sequence of $R$-modules. If $M$ is  $S$-flat, then $K$ is  $S$-flat if and only if $L$ is  $S$-flat.
\end{cor}

\proof For any $R$-module $N$, there exists an exact sequence
 $$\Tor_R^{2}(M, N_S)\to\Tor_R^{1}(K, N_S) \to\Tor_R^1 (L, N_{S}) \to\Tor_R^1(M,N_{S}).$$
Since $M$ is $S$-flat, the flanking terms are {0}, so that $\Tor_R^1 (K, N_S)\cong\Tor_R^1 (L, N_{S}).$
Therefore, by Proposition \ref{2-prop-S-flat-tor-chara}, if one of the modules $K$ and $L$ is $S$-flat, then so is the other. \cqfd

 Recall that an $R$-module $M$ is said to be pure injective provided that the induced sequence $0 \to \Hom_R(C, M) \to \Hom_R(B, M) \to \Hom_R(A, M) \to 0$ is exact for any pure exact sequence $0 \to A \to  B \to C \to  0$. We also say that $M$ is injective with respect to pure exact sequences. In this paper, we are interested in the injectivity of $M$ with respect to an other class of exact sequences. This is why we introduce the following notions:
\begin{defn}\label{2-def-S-pure-inj}
\begin{enumerate}
\item A short exact sequence of $R$-modules $0\to A\to B\to C\to 0$ is said to be  $S$-pure if the induced sequence $0\to A_S\to B_S\to C_S\to 0$ is a pure exact sequence of $R_{S}$-modules.
\item An $R$-module $M$ is said to be $S$-pure injective if it is injective relative to $S$-pure short exact sequences.
\end{enumerate} 
\end{defn}
 
\begin{rems}\label{2-rem-S-pure-inje}
\begin{enumerate} \item An $R$-module $M$ is  $S$-flat if and only if every exact sequence of $R$-modules ending with $M$ is  $S$-pure.
\item Every pure exact sequence is  $S$-pure. The converse does not hold. Indeed, le $M$ be an  $S$-flat module which is not flat \cite{Qi1}. Since $M$ is not flat there exists an exact sequence of $R$-modules $\E$, ending with $M$, which is not pure. However, due to the flatness of $M_{S}$ as an $R_{S}$-module, $\E$ is $S$-pure.
\item Every injective $R$-module is  $S$-pure injective and every  $S$-pure injective $R$-module is pure injective.
\end{enumerate}
\end{rems}


\begin{prop}\label{2-prop-S-flat-in-s-pure-seq} Let $0\to A\to B\to C\to 0$ be an $S$-pure exact sequence of $R$-modules. If $B$ is $S$-flat, then so is $C$.
\end{prop}

\proof  This follows by \cite[Theorem 1.2.14]{Gla1} and \cite[Proposition 2.6]{Qi1}. \cqfd

 \begin{prop}\label{2-prop-thecharactre-of-Rs-is-S-pure-inj} $\Hom_\Z(N,\Q/\Z)$ is  $S$-pure injective for any $R_S$-module N. \end{prop}
 
\proof Let $\E$ be an  $S$-pure exact sequence and $N$ an $R_S$-module. The result follows from the natural isomorphisms: $$\Hom_R(\E,\Hom_\Z(N,\Q/\Z))\cong \Hom_\Z(\E\otimes_R N,\Q/\Z)$$ 
 $$\E\otimes_R N\cong \E\otimes_R N_S \cong\E_S\otimes_{R_S} N_S$$
 and the fact that  $\Hom_\Z(\E\otimes_R N,\Q/\Z)$ is exact if and only if $\E\otimes_R N$ is exact \cite[Lemma 3.53]{Rot1}.\cqfd
 
\begin{cor}\label{2-cor-thecharactre-of-Rs-is-S-pure-inj} 
 $N_S^{++}:= \Hom_\Z(N_S,\Q/\Z)^+$ and $N_S^{+++}:=(N_S^{++})^+$ are  $S$-pure injective for any $R$-module $N$.
\end{cor}
\proof Obvious.\cqfd

\section{S-cotorsion modules}
 In this section, we introduce and investigate the concept of $S$-cotorsion modules. 

\begin{defn}\label{3-def-S-cotorsion}
An $R$-module $M$ is called  $S$-cotorsion if $\Ext^{1}_{R}(F,M) =0$ for any  $S$-flat $R$-module $F$. We denote by $S\C$ the class of all $S$-cotorsion modules. 
\end{defn}

\begin{rems}\label{3-rems-S-cotorsion}
\begin{enumerate}
\item Every $S$-cotorsion $R$-module is cotorsion.
\item Every injective $R$-module is $S$-cotorsion.
\item Every $S$-pure injective $R$-module is $S$-cotorsion.
\end{enumerate}
\end{rems}
\proof 1. and 2. are obvious. 

3. Let $M$ be an $S$-pure injective $R$-module and $F$ be an $S$-flat $R$-module. Let $$\E:\, 0\to K\to P\to F \to 0$$ be a short exact sequence with a projective $R$-module $P$. Consider the induced exact sequence: 
$$0\to \Hom_R(F,M)\to\Hom_R(P,M)\to\Hom_R(K,M) \to  \Ext^1_R(F,M)\to 0$$
Since $F$ is $S$-flat, $\E$ is $S$-pure. Hence,  the homomorphism $$\Hom_R(P,M)\to\Hom_R(K,M)$$ is an epimorphism because $M$ is $S$-pure injective. Hence, $\Ext^1_R(F,M)=0$. \cqfd

\begin{prop}\label{3-prop-S-cotorsion-properties}
$S\C$ is closed under extensions, direct summands and direct products.
\end{prop}
 \proof The closedness under extensions is given by the long exact sequence.  For the closedness with respect direct summands and direct products,  we use the natural isomorphism $\Ext_R^{1}(M,\prod\limits_{i\in I} C_{i})\cong\prod\limits_{i\in I}\Ext_R^{1}(M, C_{i})$ \cite[Proposition 7.22]{Rot1}. \cqfd
 
 The next result shows that, as in the case of cotorsion modules, the class of all $S$-cotorsion modules is closed under the cokernel of monomorphisms.
\begin{prop}\label{3-prop-S-cotorsion-quotion-closed}
In an exact sequence of $R$-modules, $0\to A \to B \to C\to 0$, if  both $A$ and $B$ are $S$-cotorsion, then so is $C$.
\end{prop}

\proof We claim that $\Ext^n_R(F,C') = 0$ for all $n \geq 2$ if $C'$ is $S$-cotorsion and $F$ is $S$-flat. Take a partial projective resolution of $F$ $$0 \to K \to P_{n-2} \to ... \to P_0 \to F \to 0.$$
Since $K$ is $S$-flat, $\Ext^1(K, C') = 0$. It follows from \cite[Proposition 8.5]{Rot1} that $\Ext^n_R(F, C') \cong \Ext^1_R(K, C') = 0$. 

 Now the result can be easily proved by applying the long exact sequence associated with $0\to A \to B \to C\to0$.   \cqfd 


It is well-known that all modules have a cotorsion envelope and a flat cover. The corresponding results are also true if we consider $S$-cotorsion and $S$-flat modules.

\begin{prop}\label{3-prop-S-flat-cover-S-cotorsion-envelope}
Every $R$-module has $S\F$-covers and $S\C$-envelopes. In particular,
all $R$-modules have special $S\F$-precovers and special $S\C$-preenvelopes.
\end{prop}
\proof Every $R$-modules have $S\F$-covers and $S\C$-envelopes by Lemma \ref{2-lem-S-flat-class-properties}, Lemma  \ref{3-prop-S-cotorsion-properties} (3) and \cite[Theorem 3.4]{Hol1}. The rest follows by Wakamatsu's Lemmas \cite[2.1]{Xu1}. \cqfd
 
\begin{lem}\label{3-lem-Hom-of-R-R'-is-S-cotorsion} Let $R'$ be a commutative ring. If  $E$ is an $R$-$R'$-bimodule and injective as an $R'$-module, then $\Hom_{R'}(M, E)$ is an $S$-cotorsion $R$-module for every $R_{S}$-module $M$.
\end{lem}

\proof Let $F$ be an $S$-flat $R$-module. Since $E$ is injective, there is an isomorphism \cite[Theorem 1.1.8]{Gla1}
$$\Ext^{1}_{R}(F, \Hom_{R'}(M, E))\cong  \Hom_{R'}(\Tor_R ^1(F, M), E), $$
the right hand is zero by  Proposition \ref{2-prop-S-flat-tor-chara}. \cqfd

\begin{cor}{\label{3-cor-the-chara-of-Rs-module-S-cotorsion}} 
$\Hom_{\Z}(M, \Q/\Z)$ is an $S$-cotorsion $R$-module for every $R_{S}$-module $M$.
\end{cor}

\begin{thm}\label{3-thm-(SF,SC)-perfect-cotorsion-pair}
$(S\F,S\C)$ is a hereditary perfect cotorsion pair.
\end{thm}
 \proof
 By  Proposition \ref{3-prop-S-cotorsion-properties}, to prove that $(S\F,S\C)$ is a cotorsion pair, it suffices to show that if $\Ext^1_R(F, C) = 0$ for any $S$-cotorsion $C$, then $F$ is $S$-flat.  But this means that $\Ext^1_R(F, \Hom_{\Z}(M_S, \Q/\Z)) = 0$ for every $R$-module $M$    (by Corollary \ref{3-cor-the-chara-of-Rs-module-S-cotorsion}).  Since $\Q/\Z$ is injective as an abelian group, there is an isomorphism: 
$$\Ext^{1}_{R}(F, \Hom_{\Z}(M_S, \Q/\Z))\cong  \Hom_{\Z}(\Tor_R ^1(F, M_S), \Q/\Z)$$
\cite[Theorem 1.1.8]{Gla1}. It follows that $\Tor_R ^1(F, M_S) = 0$ for every $R$-module $M.$ Therefore, by Proposition \ref{2-prop-S-flat-tor-chara}, $F$ is $S$-flat, as desired.

Finally, $(S\F, S\C)$ is hereditary by Proposition \ref{3-prop-S-cotorsion-quotion-closed} and perfect by Proposition \ref{3-prop-S-flat-cover-S-cotorsion-envelope}. \cqfd

As for the class of cotorsion modules, the class of $S$-cotorsion modules is useful when characterizing rings. Recall from \cite[Definition 4.1]{Bou1} that a ring $R$ is said to be $S$-perfect if every $S$-flat $R$-module is  projective. The following result can be viewed as  an $S$-version of \cite[Proposition 3.3.1]{Xu1}.

\begin{prop}\label{3-prop-S-perfect} Let $R$ be a commutative ring and $S$ a multiplicative subset of R. Then the following are equivalent:
\begin{enumerate}
\item $R$ is $S$-perfect.
\item Every $R$-module is $S$-cotorsion.
\item Every $S$-flat $R$-module is $S$-cotorsion.
\end{enumerate}
\end{prop}

\proof $1. \Rightarrow 2.$ and $2. \Rightarrow 3.$ are trival.\\
 $3. \Rightarrow 1.$ For any $S$-flat $R$-module $F$, we have an exact sequence $0 \to K \to P\to F\to 0$ with $P$ projective and $K$  $S$-flat. By 2. $K$ is $S$-cotorsion, and then this sequence is split. This means that $F$ is projective.  \cqfd

 The following  proposition, which may be viewed as the dual of the previous Proposition \ref{3-prop-S-perfect}, will be needed later.  

\begin{prop}\label{3-prop-dual-S-perfect} The following assertions are equivalent:
\begin{enumerate}
\item Every $R$-module is $S$-flat;
\item Every $S$-cotorsion $R$-module is $S$-flat;
\item Every $S$-pure injective $R$-module is $S$-flat.
\end{enumerate}
\end{prop}
 
\proof $1. \Rightarrow 2.$ and $2. \Rightarrow 3.$ are trivial.\\
 $3. \Rightarrow 1.$ For a fixed $R$-module $N$, we denoted by $N_S^+$, $N_S^{++}$ and $N_S^{+++} $ the $R_S$-modules $\Hom_\Z(N_S,\Q/\Z)$, $\Hom_\Z(N_S^+,\Q/\Z)$ and $ \Hom_\Z(N_S^{++},\Q/\Z)$, respectively. Let $M$ be an $R$-module. Since $\Q/\Z$ is injective as an abelian group, there is an isomorphism \cite[Theorem 1.1.8]{Gla1}
$$\Ext^{1}_{R}(M_S, N_S^{+++})\cong  ( \Tor_R ^1(M_S, N_S^{++}))^+.$$

By Corollary \ref{2-cor-thecharactre-of-Rs-is-S-pure-inj} $N_S^{++}$ is $S$-pure injective, so it is $S$-flat by 3. Hence the right hand in the above isomorphism is zero. On the other hand, $N_S^+$ is pure injective by Remark \ref{2-rem-S-pure-inje}, then it is a direct summand of $N_S^{+++}$ \cite[Proposition 2.3.5]{Xu1}. Thus $\Ext_R^1(M_S,N_S^+)=0$, then by the natural isomorphism $$\Ext^{1}_{R}(M_S,N_S^+)\cong\Tor_R ^1(M_S, N_S)^+,$$ $ \Tor_R ^1(M_S, N_S)=0$. This shows that, by Proposition \ref{2-prop-S-flat-tor-chara}, $N_S$ is $S$-flat. By Lemma  \ref{2-lem-localiz-isomorphism},  $(N_S)_S\cong  N_S$ is a flat $R_S$-module; hence $N$ is $S$-flat, as desired. \cqfd

Recall that a cotorsion envelope $\sigma_M : M \to C(M)$ has the unique mapping property \cite{Din1} if, for any homomorphism $f : M \to N$ with $N$ cotorsion, there exists a unique $g : C(M) \to N$ such that $g\sigma_M = f$. The concept of $S$-cotorsion envelopes with the unique mapping property can be defined similarly. We have the following result which may be viewed as the $S$-counterpart of \cite[Theorem 2.18]{Mao1}.

\begin{thm}\label{3-thm-S-perfect} The following assertions are equivalent.
\begin{enumerate}
\item $R$ is $S$-perfect.
\item Every $R$-module has an $S$-cotorsion envelope with the unique mapping property.
\item Every $S$-flat $R$-module has an $S$-cotorsion envelope with the unique mapping property.
\item For any $R$-homomorphism $ f : M \to N$ with $M$ and $N$ (injective) $S$-cotorsion, $ker(f)$ is $S$-cotorsion.
\item For each ($S$-flat) $R$-module $M$, the functor $\Hom_{R}(-, M)$ is exact with respect
to each $S$-pure exact sequence $0 \to K \to P \to L \to 0$ with $P$ projective.
\end{enumerate}
\end{thm}

\proof 
By using the Proposition \ref{3-prop-S-perfect} instead of \cite[ Proposition 3.3.1 ]{Xu1}, the proof is similar to that of \cite[Theorem 2.18]{Mao1}.  However, for the sake of completeness we give its proof here.

$1.\Rightarrow 2.$ and $2.\Rightarrow 3. $   are obvious.
  
  $1.\Rightarrow 4.$ This follows from Proposition \ref{3-prop-S-perfect}.
  
 $3.\Rightarrow 1.$ Let $M$ be any $S$-flat $R$-module. There is the  commutative diagram with exact arrow:
$$\xymatrix @!0 @R=20mm  @C=1.5cm{0 \ar[r]&  M \ar[drrrr]_{0}^{}  \ar[rr]^{\sigma_M}& &  S\C(M)  \ar[drr]^{\sigma_L\circ \gamma} \ar[rr]^{\gamma} & & L \ar[r]^{} \ar[d]^{\sigma_L} &  0  \\
  & &  &   &     &  S\C(L)}$$
  
Note that $\sigma_L \gamma\sigma_M = 0 = 0\sigma_M$, so $\sigma_L\gamma = 0$ by (3). Therefore $L = \Im(\gamma)\subseteq \ker(\sigma_L) = 0$, and so $M$ is  $S$-cotorsion. Hence 1. follows by Proposition \ref{3-prop-S-perfect}.

$4.\Rightarrow 1. $ Let $M$ be any $S$-flat $R$-module. We have an exact sequence $0\to M\to E\to F$ with $E$ and $F$ injective. Then, by (4), $M$ is $S$-flat.

$1.\Rightarrow 5.$ This is clear since, by Proposition \ref{2-prop-S-flat-in-s-pure-seq},  $L$  is $S$-flat.

$5.\Rightarrow 1.$ Let $M$ and $N$  be $S$-flat $R$-modules. There exists an exact sequence $0 \to K \to P\to N \to  0$ with $P$ projective, which induces an exact sequence
$$\Hom(P,M) \to \Hom_R(K, M)\to  \Ext^1(N, M)\to 0 \,\,\, (*)$$
Since $N$ is $S$-flat, $0 \to K \to P\to N \to  0$ is $S$-pure by Remarks \ref{2-rem-S-pure-inje}. Hence, by (5), $\Hom(P, M)\to  \Hom(K, M)$ is an epimorphism; so, by $(*)$, $\Ext^1(N,M)= 0$. Thus, $M$ is $S$-cotorsion. Then, (1) follows by Proposition \ref{3-prop-S-perfect}.  \cqfd
  \cqfd  

Next, we  prove the following theorem which characterizes when $R_S$ is a von Neumann regular ring. 
 
\begin{thm}\label{3-thm-dual-S-perfect} The following assertions are equivalent:

\begin{enumerate}
\item $R_{S}$ is von Neumann regular.
\item Every $S$-cotorsion $R$-module is injective.
\item Every $R$-module has an $S$-flat cover with the unique mapping property.
\item Every $S$-cotorsion $R$-module has an $S$-flat cover with the unique mapping property.
\item For any $R$-homorphism $ f : M \to N$ with $M$and $N$ (projective) $S$-flat, \coker(f) is $S$-flat.
\end{enumerate}
\end{thm}

\proof Notice that $R_S$ is von Neumann regular if and only if every $R$-module is $S$-flat.

 $1.\Leftrightarrow 2.$ This follows from Definition \ref{3-def-S-cotorsion} and Proposition \ref{3-prop-dual-S-perfect}.

$1.\Rightarrow 3.,$ $3.\Rightarrow 4.$ and $1. \Rightarrow 5.$ are obvious.

 $4.\Rightarrow 1.$  
Let $M$ be any $S$-cotorsion $R$-module. There is the exact commutative diagram.
$$\xymatrix @!0 @R=15mm  @C=2.5cm {  &   S\F(K) \ar[drr]^{0} \ar[d]_{\varepsilon_K} \ar[dr]_{\alpha\circ \varepsilon_K} &  &  \\
 0 \ar[r]{} &K \ar[r]_{\alpha} \ar[d]{}& S\F(M)\ar[r]_{\varepsilon_M}  &  M  \ar[r]{}   & 0    \\
  &  0 &   &     &   }$$
Note that $\varepsilon_M\alpha\varepsilon_K = 0 = \varepsilon_M 0$, so $\alpha \varepsilon_K = 0$ by $(4)$. Therefore
 $K = \Im(\varepsilon_K)\subseteq \ker(\alpha) =0$, and so $M$ is $S$-flat. Hence, $R_S$ is von Neumann regular by Proposition \ref{3-prop-dual-S-perfect}.  
 
  $5.\Rightarrow 1.$ Let $M$ be any $R$-module. We have an exact sequence $P\to Q\to M\to 0$ with $P$ and $Q$ projective. Then, by (5), $M$ is $S$-flat. \cqfd  

Next, following \cite{Mao1}, we focus our attention on the question when $N$ and $M$ share a common $S$-cotorsion (pre)envelope. 

   
 
 
\begin{prop}\label{3-prop-S-envelope} Let $\alpha : N \to M $ be a monomorphism. Then we have:

 \begin{enumerate}
\item If $ \coker(\alpha)$ is $S$-flat, then $\beta\alpha : N \to H$ is an $S$-cotorsion preenvelop of $N$ whenever $\beta : M \to  H$ is an $S$-cotorsion preenvelope of $M$.
\item $\sigma_{M}\alpha : N \to SC(M)$ is a special $S$-cotorsion preenvelope of $N$ if and only if $ \coker(\alpha)$ is $S$-flat. 
 \end{enumerate}
 
\end{prop}
\proof The proof is similar to that of \cite[ Proposition 2.6]{Mao1}.\cqfd


\begin{prop} Assume the class of $S$-flat $R$-modules is closed under cokernels of monomorphisms. If $\alpha : N \to M$ is an essential monomorphism with $ \coker(\alpha)$ $S$-flat, then there exists $h : M \to S\C(N)$ such that $h$ is a special $S$-cotorsion preenvelope of $M$.
\end{prop}
 
\proof The proof is similar to that of \cite[ Proposition 2.7]{Mao1}. \cqfd
 As is  well-known, for two right $R$-modules $N\leq_e M$, if $M \leq_e \C(N)$, then $\C(N)= \C(M)$ if and only if $M/N$ is flat if and only if $\C(M)/N$ is flat \cite[Theorem 2.8]{Mao1}. Replacing "cotorsion" with "$S$-cotorsion", we have
\newpage

\begin{thm}
Assume that $N\leq_{e}  M \leq_{e} SC(M)$. Then the following are equivalent:
\begin{enumerate}
\item $M/N$ is $S$-flat;
\item $S\C(M)/N$ is $S$-flat;
\item $S\C(M)= S\C(N)$ (up to isomorphism).
\end{enumerate}
\end{thm}
\proof We imitate the proof given by \cite[Theorem 2.8]{Mao1} with some changes. 

$1. \Rightarrow 3.$ Let $i : N\to M$ be the inclusion map. Since $M/N$ is $S$-flat,
there is  $\alpha : M\to S\C(N)$ such that $\alpha i = \sigma_N$. Note that $\alpha$ is monic since $\sigma_N$ is monic
and $i$ is an essential monomorphism \cite[Corollary 5.13]{And3}. By the defining property of an $S$-cotorsion
envelope, it follows that $\alpha$ factors through $\sigma_M: M \to S\C(M)$, so there is $f : S\C(M)\to S\C(N)$ such that $f\sigma_M = \alpha$. Since $\alpha$ is monic and $\sigma_M$ is an essential monomorphism, $f$ is a monomorphism. Similarly, the map $\sigma_M i : N \to S\C(M)$ factors through $\sigma_N : N \to S\C(N)$, so
there is $g : S\C(N) \to S\C(M)$ such that $\sigma_M i = g\sigma_N$. Thus $\sigma_N = \alpha i = f\sigma_M i  = fg\sigma_N$, which implies $fg$ is an automorphism of $S\C(N)$ by the defining property of an $S$-cotorsion
envelope, and hence $f$ is an epimorphism. It follows that $f$ is an isomorphism.

$3. \Rightarrow 2.$ This is clear.

$2. \Rightarrow 1.$ There is an exact sequence

$$0 \to M/N\to S\C(M)/N\to S\C(M)/M\to 0.$$
The $S$-flatness of $S\C(M)/N$ and $S\C(M)/M$ implies that $M/N$ is $S$-flat by Lemma \ref{2-lem-S-flat-class-properties}.
\cqfd


In \cite[Proposition 2.9]{Mao1}, the authors show that if $N \leq M \leq \C(N)$ and $M \leq_{e} \C(M)$, then $\C(M) = \C(N)$ (up to isomorphism), where $C(M)$, $\C(N)$ are the cotorsion envelopes of $N$ and $M$, respectively. Here we have the corresponding result for $S$-cotorsion envelopes.

\begin{prop} 
If $N \leq  M \leq S\C(N)$ and $M \leq_{e} S\C(M)$, then $S\C(M) = S\C(N)$ (up to isomorphism).
\end{prop}

\proof The proof is similar to that of \cite[ Proposition 2.9]{Mao1}
\section{$S$-flat  preenvelopes and $S$-coherent rings}
Recall from \cite[Definition 3.3]{Ben1} that a ring $R$ is called $S$-coherent, where $S$ is a multiplicative subset of $R$, if every finitely generated ideal of $R$ is $S$-finitely presented. In this section, we demonstrate a new characterization of these rings in terms of $S\F$-preenvelopes. We cite this lemma here:

\begin{lem}[\cite{Xu1}, Lemma  2.5.2]\label{4-lem}
For any ring $R$, if $N\subseteq M$ be a submodule, then $N$ can be enlarged to a
submodule $N^*$ such that $N^*$ is pure in $M$ and the cardinality of $N^*$ is less than or equal to $ \Card(N) \Card(R)$ if either of $ \Card(N)$ and $ \Card(R)$ is infinite. If both are finite, there is an $N^*$ which is at most countable.
\end{lem}

\begin{thm}\label{4-thm-S-flat-product-envelope} $S\F$ is closed under direct products if and only if every $R$-module has an $S\F$-preenvelopes.
\end{thm}
\bigskip

\proof We imitate the proof given by \cite[Theorem 2.5.1]{Xu1} with some changes. 

$1.\Rightarrow 2.$ For any  $R$-module $M$, let $ \Card(M) \Card(R)\leq \aleph_\beta$, where $\aleph_\beta$ is an infinite cardinal number. Set $$\X = \lbrace G\in S\F| \Card(G)\leq \aleph_\beta\rbrace.$$
  Let $(G_i)_{i \in I}$ be a family of representatives of this class with the index set $I$. Let $H_i = \Hom_R(M, G_i)$ for each $i\in I$ and let $F =\prod G_i^{H_i}$. Define $\varphi : M\to F$ so that the composition of $\varphi$ with the projection morphism $F\to G_i^{H_i}$ maps $x\in F$ to $(h(x))_{h\in H_i}$. By assumption $F$ is an $S$-flat $R$-module. We claim that $\varphi : M \to F$ is an $S$-flat preenvelope of $M$. Let $\varphi' : M \to G$ be a linear map with $G$ $S$-flat. By  Lemma \ref{4-lem}, the submodule $\varphi'(M)\subseteq G$ can be enlarged to a pure submodule $G'\subseteq G$ with $ \Card(G')\leq \aleph_\beta$. Notice that $G'$ is $S$-flat because it is a pure submodule of $G$ which is $S$-flat. Then $G'$ is isomorphic to one of the $G_i$. By the construction of the morphism $\varphi$, $\varphi'$ factors through $\varphi$, as desired.

$2.\Rightarrow 1.$ 
Let $(F_i)_{i\in I}$ be a family of $S$-flat modules, and $\prod\limits_{i\in I} F_i\to F$ be an  $S\F$-preenvelope. Then there are factorizations $\prod\limits_{i\in I} F_i\to F\to F_i$ (induced by the canonical projection $\prod\limits_{i\in I} F_i\to F_i$). These give rise to a map $F\to \prod\limits_{i\in I} F_i$  with the composition $\prod\limits_{i\in I} F_i\to F \to \prod\limits_{i\in I} F_i$  the identity. Hence $ \prod\limits_{i\in I} F_i$ is isomorphic to a summand of $F$ and so it is $S$-flat.  \cqfd

Recall from \cite[Theorem 4.4]{Qi1} that a ring $R$ is $S$-coherent if and only if the direct product of any family of flat $R$-modules is $S$-flat. We have the following consequence.
\begin{cor}\label{4-cor-S-coherent}
Assume that every module has an $S\F$-preenvelope, then $R$ is $S$-coherent. 
\end{cor}

\proof Using the Theorem  \ref{4-thm-S-flat-product-envelope} and the fact that any flat module is $S$-flat, the result follows from \cite[Theorem 4.4]{Qi1}.\cqfd 

It is worth noting that in \cite[Theorem 2.5.1]{Xu1} coherent rings are characterized
using the notion of flat preenvelopes. Namely, for a ring $R$, every left $R$-module $M$ has a flat preenvelope if and only if $R$ is right coherent. Naturally, one can ask for an $S$-version of this result, representing the converse of Corollary \ref{4-cor-S-coherent}. We leave this as an interesting open question, and here we provide a partial response.

\begin{prop}\label{4-prop-S-coherent-Rs-EP}
Assume that $R_S$ is finitely presented as an $R$-module. If $R$ is $S$-coherent, then every module has an $S\F$-preenvelope.
\end{prop}

\proof By Theorem \ref{4-thm-S-flat-product-envelope}, it suffices to show that $S\F$ is  closed under direct products. Let $(F_i)_{i\in I}$ be a family of $S$-flat modules. We need to prove that $(\prod F_i)_S=R_S\otimes_R\prod F_i$ is a falt $R_S$-module. Since $R_S$ is finitely presented, then $R_S\otimes_R\prod F_i \cong \prod R_S\otimes_R F_i$ as $R$-modules \cite[Theorem 3.2.22]{Eno2}. Note that $R_S\otimes_R F_i$ is a falt  $R$-module for any $i\in I $. Since $R$ is $S$-coherent, then $ \prod R_S\otimes F_i$ is $S$-flat \cite[Theorem 4.4]{Qi1}. Hence, $(\prod F_i)_S$ is a flat $R_S$-module and so $\prod F_i$ is $S$-flat, as desired. \cqfd 
 We end this paper with the following consequence, which shows that the concepts of $S$-coherent and coherent rings coincide on $S$-perfect rings.
\begin{cor}
Let $R$ be an $S$-perfect ring. Then $R$ is $S$-coherent if and only if it is coherent.
\end{cor}
\proof Recall from \cite[Proposition 3.5]{May1} that a ring $R$ is coherent and perfect if and only if every $R$-module has a projective (pre)envelope. 

The "if"  part always true. 

The "only if"  part  follows form Proposition \ref{4-prop-S-coherent-Rs-EP} and the fact that $R_S$ is a finitely presented $R$-module whenever $R$ is an $S$-perfect ring  \cite[Theorem 3.10]{Bou1}. \cqfd

\begin{exmp}
Let $R_1$ be an $S_1$-perfect coherent ring (semisimple ring as an example), $R_2$ be any commutative ring which is not coherent. Consider the ring $R=R_1\times R_2$ with the multiplicative subset $S = S_1\times 0$. Then 

\begin{enumerate}
\item $R_S\cong (R_1)_{S_1}\times 0$ is a finitely presented projective $R$-module. 
\item $R$ is an $S$-coherent ring which is not coherent.
\end{enumerate}  
\end{exmp}

\proof
1. Since $R_1$ is $S_1$-perfect, $(R_1)_{S_1}$ is a finitely generated projective $R_1$-module by \cite[Theorem 4.9]{Bou1}. Then  $R_S\cong (R_1)_{S_1}\times 0$ is a finitely generated projective $R$-module; so, it is finitely presented. 

2. $R$  is $S$-coherent by \cite[Proposition 3.5]{Ben1}. \cqfd
\noindent\textbf{Acknowledgment.} The authors wish to express their gratitude to the referee for the careful critical reading of the manuscript and for his/her valuable comments.

Driss Bennis:  Faculty of Sciences, Mohammed V University in Rabat, Rabat, Morocco.

\noindent e-mail address: driss.bennis@um5.ac.ma; driss$\_$bennis@hotmail.com

Ayoub Bouziri: Faculty of Sciences, Mohammed V University in Rabat, Rabat, Morocco.

\noindent e-mail address: ayoub$\_$bouziri@um5.ac.ma

\end{document}